\documentclass{article}[12pt]
\usepackage{myfig,latexsym,amsthm,amsfonts,epsfig,psfrag
}
\begin{document}

\newtheorem{lemma}{Lemma}
\newtheorem{theorem}{Theorem}

\def \H{{I\!\!H}}
\def \R{{I\!\!R}}
\def \E{{I\!\!E}}
\def \Z{{\mathbb Z}}

\begin{center}
{\Large  \bf
Hyperbolic Coxeter $n$-polytopes\\ with $n+2$ facets
}

\vspace{15pt}
{\large P.~Tumarkin}
\end{center}

\vspace{8pt}

\begin{center}
\parbox{10cm}%
{\scriptsize 
{\it Abstract.}
In this paper, we classify all the hyperbolic non-compact 
Coxeter polytopes of finite volume combinatorial type of
which is either a pyramid over a product of two simplices 
or a product of two simplices of
dimension greater than one. Combined with results of
Kaplinskaja~\cite{Kap} and Esselmann~\cite{Ess} this
completes the classification of hyperbolic Coxeter $n$-polytopes
of finite volume with $n+2$ facets. 
}
\end{center}

\vspace{8pt}

\section{Introduction}

Consider a convex polytope $P$ in the $n$-dimensional hyperbolic space
 $\H^n$. 

A polytope is called a {\it Coxeter polytope} if its dihedral angles are
 integer parts of $\pi$. Any Coxeter polytope is a fundamental domain
of the discrete group generated by the reflections with respect to
its facets.

We are interested in Coxeter polytopes of finite volume. In fact,
there is no complete classification of such polytopes, but some
partial results are obtained. For example, Coxeter polytopes in
$\H^3$ are completely characterized by Andreev
\cite{A1},~\cite{A2}; all the hyperbolic simplices are classified
 (\cite{L},~\cite{V3}).

There are some results about hyperbolic $n$-polytopes with $n+2$
facets. Kaplinskaja \cite{Kap} described all Coxeter simplicial
prisms, and Esselmann \cite{Ess} obtained the classification of
compact Coxeter hyperbolic $n$-polytopes with $n+2$ facets.
In this paper, we complete the classification of Coxeter 
hyperbolic $n$-polytopes of finite volume with $n+2$ facets.

The author is grateful to Prof. R.~Kellerhals who brought 
the work~\cite{Ess} to his attention.  

\section{Gale diagrams}
\label{sec-Gale}

The essential facts about convex polytopes and Gale diagrams 
which are used in this paper
are quoted below. Proofs and details may be found in~\cite{gale}. 

Every combinatorial type of $n$-polytope with $k$ facets can be
represented by a Gale diagram. This consists of $m$ points 
$a_1,\dots,a_m$ of the unit sphere in $\R^{k-(n+1)}$ centered at the
origin and (possibly) the 
origin $a_0$ which are labeled according to the following rules:

1) Each label is a positive integer, and the sum of labels equals $k$.

2) The points that lie in any open half-space bounded by a hyperplane
 through the origin have labels whose sum is at least two.

The combinatorial type of a convex polytope can be read off from 
the Gale diagram in the following way.
Each point $a_i$ with label $\mu_i$ corresponds to $\mu_i$ facets 
$f_{i,1},\dots,f_{i,\mu_i}$ of $P$. For any subset $I$ of the set of
facets of $P$ the intersection of facets 
$\{f_{j,\gamma} | (j,\gamma)\in I \}$ 
is a face of $P$ if and only if $0$ is contained in the relative interior of 
$conv\{a_{j} | (j,\gamma)\notin I\}$.

By {\it pyramid} in $\H^n$ we mean a convex hull of a point (apex)
and an $(n-1)$-dimensional polytope which is not a simplex.
A polytope $P$ is called a {\it multipyramid } if 
$P$ is a pyramid over
a pyramid.

Now let $P$ be a $n$-polytope with $n+2$ facets.
In this case the Gale diagram is one-dimensional,
and each point in the diagram is equal to 
$-1$, $0$ or $+1$. The {\it multiplicity} of $-1$, $0$ and $+1$
is defined to be the label of correspondent point in the Gale diagram.
Suppose that the multiplicities are $p$, $q$ and $r$ respectively,
$p+q+r=n+2$.
There exists a polytope correspondent to a Gale diagram if and only if
$p,r\ge 2 $ holds. 

In the case $q=0$ the combinatorial type of $P$
is a product of two simplices of dimension $p-1$ and $r-1$.
The set of facets is a union of two disjoint sets $J_p$ and $J_r$
of order $p$ and $r$ respectively. 
Any vertex of $P$ is an intersection of $p-1$ facets from $J_p$
and  $r-1$ facets from $J_r$.

If $q=1$ then  the combinatorial type of $P$
is a pyramid over a product of two simplices.
If $q> 1$ then  the combinatorial type of $P$
is a multipyramid over a product of two simplices.

\begin{lemma}
\label{pyr-product}
Let $P$ be a Coxeter polytope in $\H^n$ with $n+2$ facets.
Then  the combinatorial type of $P$ is
either a product of two simplices or a pyramid over
a product of two simplices.

\end{lemma}

\begin{proof}
It is sufficient to prove that any 
pyramid over a pyramid is not a Coxeter polytope.

Suppose that $P$ is a pyramid over a pyramid $P'$.
Since $P'$ is a pyramid,  it is not simple.
Hence, the apex $A$ of the pyramid $P$ is an ideal vertex of $P$.
Consider a small  horosphere $h$ centered in $A$.
The combinatorial type of the intersection $h\cap P$  
coincides with the combinatorial type of $P'$.
At the same time, $h\cap P$  is a compact Euclidean Coxeter
$(n-1)$-polytope. But  any compact
Euclidean Coxeter polytope is simple.

\end{proof}

Any Coxeter polytope $P$ can be represented by its  Coxeter diagram 
$\Sigma(P)$.
Nodes of Coxeter diagram correspond to the facets of $P$. Two nodes
are joined by a $(m-2)$-fold edge or a $m$-labeled edge if the 
corresponding dihedral angle equals $\frac{\pi}{m}$. If the
corresponding facets are parallel the nodes are joined by a bold edge,
and if they diverge then the nodes are joined by a dotted edge
labeled by $cosh(\rho)$, where $\rho$ is the distance between the facets.

Let $\Sigma$ be a diagram of order $k$ with nodes $v_1$,...,$v_k$.
Define a  symmetric $k\times k$ matrix $G(\Sigma)$ with $g_{ii}=1$,
$i=1,...,k$, and $g_{ij}=-cos(\frac{\pi}{k})$ 
($g_{ij}=-1$ or $g_{ij}=-cosh(\rho)$)
if $v_i$ and $v_j$ are joined by $k$-labeled edge (by a bold edge or by 
a dotted edge labeled by $cosh(\rho)$ respectively).  

If $\Sigma$ is a Coxeter diagram of a polytope $P$ then 
the matrix $G(\Sigma)$ coincides with the Gram matrix of $P$.
See \cite{V1} for  details.

By the signature and the determinant of a diagram $\Sigma$ 
we mean the signature and the determinant of the matrix
$G(\Sigma)$.

\section{Products of two simplices}
\label{products}

In this section, we classify all the non-compact $n$-polytopes of finite
volume with $n+2$ facets,  combinatorial type of which
is a product of two simplices of dimension greater than one. 
For the classification of compact polytopes of this combinatorial type
see \cite{Ess}.
All simplicial prisms are listed in~\cite{Kap} (see also~\cite{V1}).
 
The result of this section is 

\begin{figure}[htb!]
\begin{center}
\psfrag{v1}{$1$}
\psfrag{v2}{$2$}
\psfrag{v3}{$3$}
\psfrag{v4}{$4$}
\psfrag{v5}{$5$}
\psfrag{v6}{$6$}
\epsfig{file=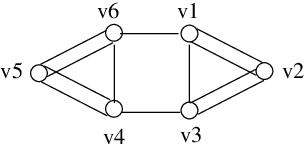,width=0.3\linewidth}
\end{center}
\caption{}
\label{prod}
\end{figure}

\begin{theorem}
There exists a unique
non-compact Coxeter polytope of finite volume 
which is combinatorially equivalent to a product of two simplices 
of dimension greater than one. It is shown in Fig.~\ref{prod}.

\end{theorem}

A diagram $\Sigma$ is called a {\it Lann\'er diagram} if any subdiagram
of $\Sigma$ is elliptic, and the diagram
$\Sigma$ is neither elliptic nor parabolic.
All Lann\'er diagrams are classified in~\cite{L}.

\begin{lemma}
\label{product}
Let $P$ be a hyperbolic  Coxeter polytope which combinatorial type 
is  a product of two
simplices and $\Sigma$ be a Coxeter diagram of $P$. 
Then $\Sigma$ satisfies the following three conditions:

\begin{itemize}
\item[{\sc (i)}]
$\Sigma$ is a union of two disjoint Lann\'er diagrams $L_1 $ and $L_2$.

\item[{\sc (ii)}]
For any $v_1\in L_1$ and  $v_2\in L_2$
a diagram $\Sigma\!\setminus\! \{v_1,v_2\}$ is either elliptic or 
parabolic.

\item[{\sc (iii)}]
$Det(\Sigma)=0$.

\end{itemize}

\end{lemma}

\begin{proof}
According to Section~\ref{sec-Gale}, $\Sigma$ is a union of two 
disjoint diagrams, say $L_1$ and $L_2$, where $L_1$ consists of
facets correspondent to $-1$ in the Gale diagram of $P$, and
$L_2$ consists of facets correspondent to $1$.
Any vertex of $P$ corresponds to a subdiagram
of $\Sigma$ that is $\Sigma\!\setminus\! \{v_1,v_2\}$, where
$v_1\in L_1$, $v_2\in L_2$. Conversely, 
for any $v_1\in L_1$ and  $v_2\in L_2$
a diagram $\Sigma\!\setminus\! \{v_1,v_2\}$ corresponds to a vertex of
$P$.
Thus, $\Sigma\!\setminus\! \{v_1,v_2\}$ is either elliptic or 
parabolic.
Any subdiagram of $\Sigma\!\setminus\! \{v_1,v_2\}$ is elliptic.
In particular, $L_1\!\setminus\! v_1$ and  $L_2\!\setminus\! v_2$
are elliptic for any $v_1\in L_1$ and $v_2\in L_2$. 

By \cite{V1}, Th.~3.1, a subdiagram of  $\Sigma$ defines a proper face
of $P$ if and only if it is elliptic. 
Hence, $L_1$ and $L_2$ are not elliptic.
By \cite{V1}, Th.~4.1,
 $L_1$ and $L_2$ are not parabolic.
Therefore,  $L_1$ and $L_2$ are Lann\'er diagrams,
and we have proved  {\sc (i)} and {\sc (ii)}.
Condition {\sc (iii)} follows immediately from the fact that
the hyperbolic $n$-polytope $P$ has more than $n+1$ facets.

\end{proof}

Conditions {\sc (i)} and  {\sc (ii)} are satisfied by a big number
of Coxeter diagrams. All the diagrams containing parabolic subdiagrams
are listed in Table~\ref{prod_bad}.
The rest diagrams can not correspond to non-compact polytope. 
To compute the determinants of diagrams shown in the left column
we use  technical tools from~\cite{V2} (Prop. 13, Prop. 15 and Table~2).
No of these polytopes has $Det=0$.

\begin{table}[htb!]
\begin{center}
\psfrag{m5}{\scriptsize $m\ge 5$}
\psfrag{m4}{\scriptsize $m\ge 4$}
\psfrag{m3}{\scriptsize $m\ge 3$}
\psfrag{l5}{\scriptsize $l\ge 5$}
\psfrag{l3}{\scriptsize $l\ge 3$}
\epsfig{file=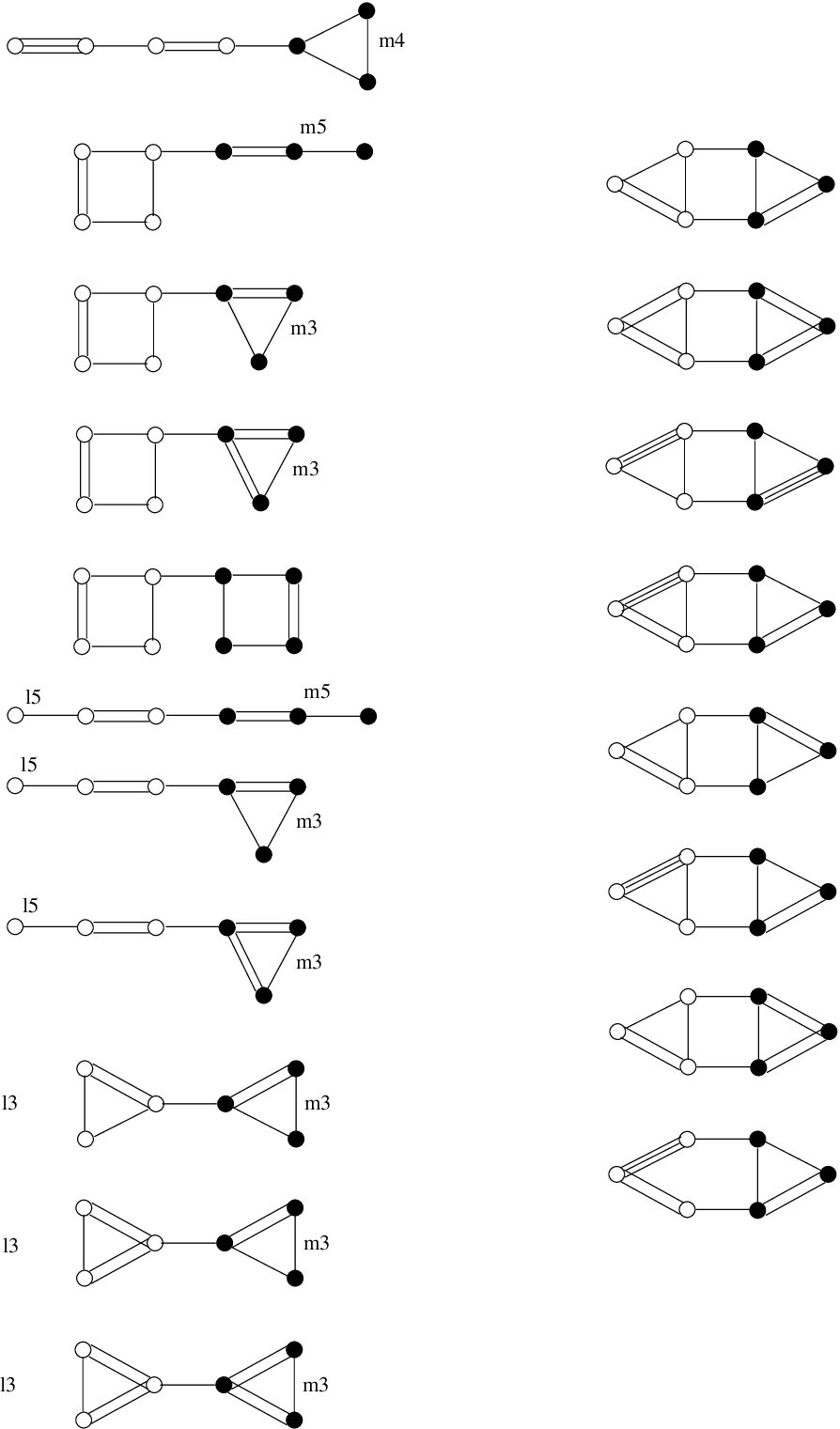,width=0.7\linewidth}
\caption{Diagrams satisfying {\sc (i)} and {\sc (ii)} and containing
parabolic subdiagrams. }
White and black vertices compose two disjoint Lann\'er subdiagrams.
\end{center}
\label{prod_bad}
\end{table}

A direct computation of other determinants shows that the only Coxeter
diagram with $Det=0$ is one shown in
Fig.~\ref{prod}.
Denote it by $\Sigma_0$. Denote its nodes as shown in Fig.~\ref{prod}.

Now we shall check the signature of
$\Sigma_0$. Since $\Sigma_0\!\setminus\! \{v_1,v_2,v_3\}$ is a Lann\'er
diagram, the negative inertia index is at least one. $\Sigma_0\!\setminus\!
\{v_1,v_4\}$ is elliptic, and the positive inertia index is 
at least 4. The determinant of  $\Sigma_0$ equals 0. Thus, the signature is
(4,1,1), where the first number is the positive inertia index, the second 
is the negative inertia index and the third one is the corank of $\Sigma_0$.

By Theorem 2.1 of~\cite{V1} for any Coxeter diagram of signature
$(n,1,l)$ there exists a polytope in $\H^n$ correspondent to this
diagram. Hence, there exists a Coxeter polytope with the Coxeter 
diagram $\Sigma_0$.  
According to~\cite{V1}, Prop.~4.2, this polytope has finite volume.

\section{Pyramids}
\label{pyramid}

In this section, we classify all the $n$-polytopes of finite
volume with $n+2$ facets,  combinatorial type of which
is a pyramid over a 
product of two simplices.

We need the following

\begin{lemma}[\cite{V2}, Prop. 12]
\label{det}
Let $\Sigma$ be a Coxeter diagram and $v$ be a node of $\Sigma$.
Suppose that $\Sigma\!\setminus\! v$ is a union of two disjoint
parabolic diagrams.
Then $Det(\Sigma)=0$.

\end{lemma}

A diagram $\Sigma$ is called a {\it quasi-Lann\'er diagram} if any subdiagram
of $\Sigma$ is either elliptic or parabolic, and the diagram
$\Sigma$ is neither elliptic nor parabolic.
The complete list of quasi-Lann\'er diagrams can be found in~\cite{V3}.

\begin{lemma}
\label{pyr}
Let $P$ be a hyperbolic  Coxeter pyramid over a product of two
simplices and $\Sigma$ be a Coxeter diagram of $P$. 
Then $\Sigma$ satisfies the following three conditions:

\begin{itemize}
\item[{\sc (i)}]
$\Sigma$ is a union of two quasi-Lann\'er diagrams $L_1 $ and $L_2$.
The intersection of $L_1$ and $L_2$ is a unique node $v$.
$L_1\!\setminus\! v$ and $L_2\!\setminus\! v$ are not adjacent.

\item[{\sc (ii)}]
The diagrams $L_1\!\setminus\! v$ and $L_2\!\setminus\! v$ are parabolic. 
Any other subdiagram of $L_1$ or $L_2$ is elliptic.

\item[{\sc (iii)}]
For any $v_1\in L_1\!\setminus\! v$ and  $v_2\in L_2\!\setminus\! v$
a diagram $\Sigma\!\setminus\! \{v_1,v_2\}$ is either elliptic or 
parabolic.

\end{itemize}

Any Coxeter diagram satisfying  {\sc (i)}--{\sc (iii)}  determines
a  hyperbolic  Coxeter pyramid over a product of two
simplices.

\end{lemma}

\begin{proof}
{\bf 1.}
Let $A$ be the apex of the pyramid $P$ and
$v$ be the node of $\Sigma$ correspondent to the facet of $P$ opposite to $A$.
Consider a small  horosphere $h$ centered in $A$.
The intersection $h\cap P$ is a direct product of two Euclidean
simplices. Thus, $\Sigma\!\setminus\! v$ is a union  of two disjoint 
parabolic diagrams, say $ S_1$ and $ S_2$.
Clearly, $ S_1$ and $ S_2$ are not adjacent.

Any vertex of $P$ except $A$ corresponds to a subdiagram 
 $\Sigma\!\setminus\! \{v_1, v_2\}$, where
$v_1\in  S_1$, $v_2\in  S_2$.
Hence,  $\Sigma\!\setminus\! \{v_1, v_2\}$ is either elliptic or
 parabolic. 
Since $ S_1$ contains at least two vertices,
 the diagram $\Sigma\!\setminus\! \{ S_1, v_2\}$ is a proper subdiagram of
 $\Sigma\!\setminus\! \{v_1, v_2\}$. Thus, it is elliptic.

Denote by $L_1$ and  $L_2$ subdiagrams
 $ S_1\!\cup\! v$ and  $ S_2\!\cup\! v$ respectively.
We have proved that for any node $u\in L_1$
the diagram $L_1\!\setminus\! u$ is either parabolic (if $u$ coincides
 with $v$) or
elliptic (otherwise).
Hence, $L_1$ is a quasi-Lann\'er diagram with the only parabolic
subdiagram $ S_1$.
Analogously,  $L_2$ is a quasi-Lann\'er diagram with the only parabolic
subdiagram $ S_2$. 

{\bf 2.}
We have shown that conditions {\sc (i)}--{\sc (iii)} are necessary.
Now, suppose that $\Sigma$ satisfies all the three conditions.
The negative inertia index of $\Sigma$ equals  1,
since $\Sigma$ has a Lann\'er subdiagram $L_1$ and the subdiagram
$\Sigma\!\setminus\! v$ is parabolic.
For any $v_1\in L_1\!\setminus\! v$
and  $v_2\in L_2\!\setminus\! v$ the subdiagram $\Sigma \!\setminus\!
\{v,v_1,v_2\}$ is elliptic. Hence, the positive inertia index of
$\Sigma$ is at least $n-1$.
By Lemma~\ref{det} $Det(\Sigma)=0$.
Thus, the signature of $\Sigma$ is either (n-1,1,2) or (n,1,1).
According to~\cite{V1}, Th. 2.1, 
there exists a Coxeter polytope in $\H^k$ with Coxeter diagram $\Sigma$,
where $k$ is either $n-1$ or $n$.

Suppose that the signature of $\Sigma$ is (n-1,1,2).
Then there exists a Coxeter polytope $P$ in $\H^{n-1}$ with Coxeter
diagram   $\Sigma$. 
This is impossible, since  
$\Sigma \!\setminus\! v$ is a parabolic subdiagram of $\Sigma$ 
of rank $n-1$. 

Hence, the signature of  $\Sigma$
equals (n,1,1), and $\Sigma$ determines a Coxeter polytope $P$ in
$\H^n$.
By Lemmas~\ref{pyr-product} and \ref{product},
 $\Sigma$ determines
a  hyperbolic  Coxeter pyramid over a product of two
simplices. 
This polytope is of finite volume by~\cite{V1}, Prop.~4.2.

\end{proof}

Considering all the quasi-Lann\'er diagrams with exactly one parabolic
subdiagram, we obtain the following 

\begin{theorem}
All the 
 Coxeter polytopes of finite volume 
the combinatorial type of which is
a pyramid over a
 product of two simplices  are listed in Tables~\ref{pyr3}--\ref{pyr17}. 
 
\end{theorem}

\noindent
{\bf Remark.} Black nodes in  Tables~\ref{pyr3}--\ref{pyr17}
correspond to the bases of pyramids.

\begin{table}[htb!]
\begin{center}
\psfrag{k}{\scriptsize $k$}
\psfrag{l}{\scriptsize $l$}
\psfrag{m}{\scriptsize $m$}
\psfrag{n}{\scriptsize $n$}
\psfrag{k=2,3,4;}{\scriptsize $k=2,3,4;$}
\psfrag{l=2,3,4;}{\scriptsize $m=2,3,4;$}
\psfrag{m=3,4;}{\scriptsize $l=3,4;$}
\psfrag{n=3,4.}{\scriptsize $n=3,4.$}
\psfrag{m=2,3.}{\scriptsize $l=2,3,4,5,6.$}
\psfrag{k=5,6;}{\scriptsize $k=5,6;$}
\psfrag{l=2,3,4,5,6;}{\scriptsize $m=2,3;$}
\epsfig{file=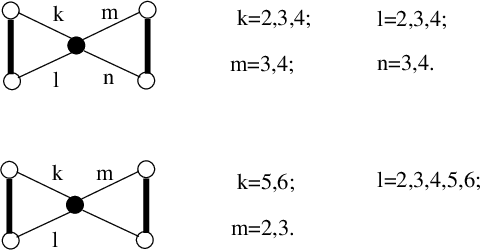,width=0.5\linewidth}
\end{center}
\caption{Pyramids in $\H^3$.}
\label{pyr3}
\end{table}

\begin{table}[htb!]
\begin{center}
\psfrag{k}{\scriptsize $k$}
\psfrag{l}{\scriptsize $l$}
\psfrag{k=2,3.}{\scriptsize $k=2,3.$}
\psfrag{k=2,3,4;l=3,4.}{\scriptsize $k=2,3,4;\ \ l=3,4.$}
\psfrag{k=2,3,4,5;l=3,4,5.}{\scriptsize $k=2,3,4,5;\ \ l=3,4,5.$}
\psfrag{k=2,3,4,5.}{\scriptsize $k=2,3,4,5.$}
\psfrag{l=3,4.}{\scriptsize  $l=3,4.$}
\psfrag{5}{\scriptsize $5$}
\epsfig{file=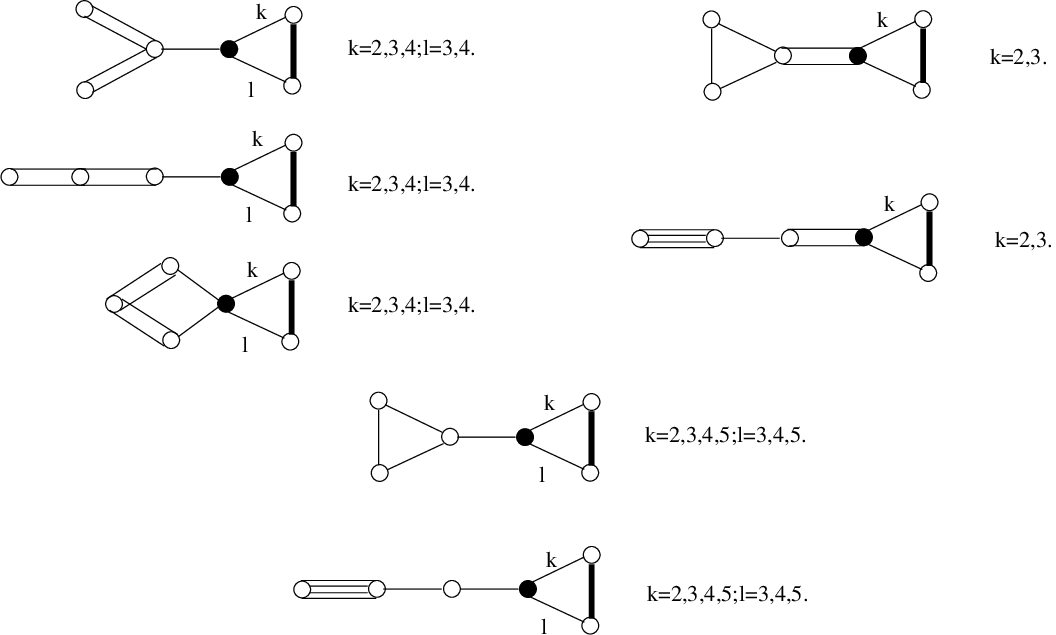,width=0.9\linewidth}
\end{center}
\caption{Pyramids in $\H^4$.}
\label{pyr4}
\vspace{40pt}
\end{table}

\vspace{40pt}

\begin{table}[htb!]
\begin{center}
\psfrag{k}{\scriptsize $k$}
\psfrag{k=2,3.}{\scriptsize $k=2,3.$}
\psfrag{k=2,3,4;l=3,4.}{\scriptsize $k=2,3,4;\ \ l=3,4.$}
\epsfig{file=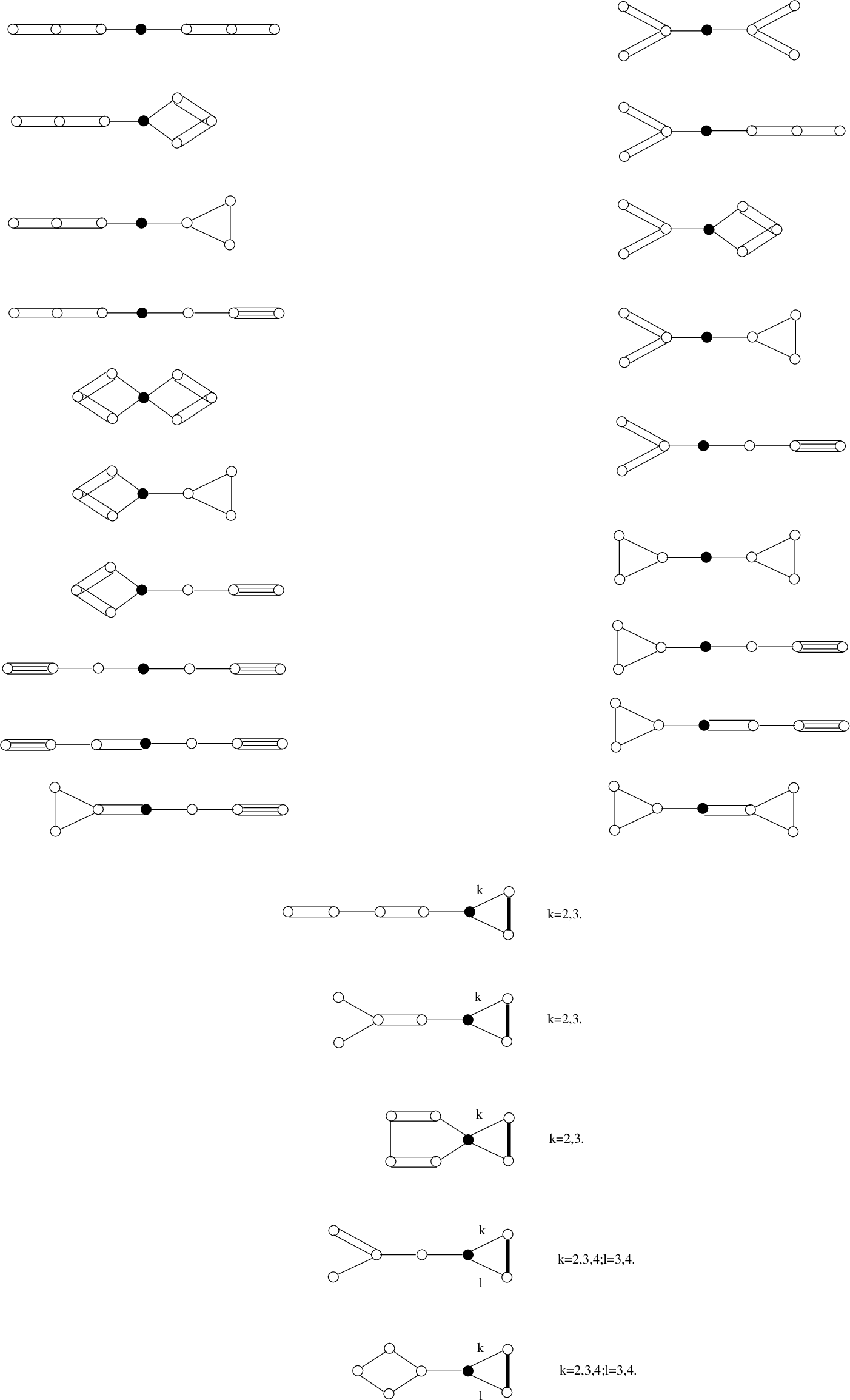,width=0.9\linewidth}
\end{center}
\caption{Pyramids in $\H^5$.}
\label{pyr5}
\end{table}

\begin{table}[htb!]
\begin{center}
\psfrag{k}{\scriptsize $k$}
\psfrag{l}{\scriptsize $l$}
\psfrag{k=2,3.}{\scriptsize $k=2,3.$}
\psfrag{k=2,3,4;l=3,4.}{\scriptsize $k=2,3,4;\ \ l=3,4.$}
\epsfig{file=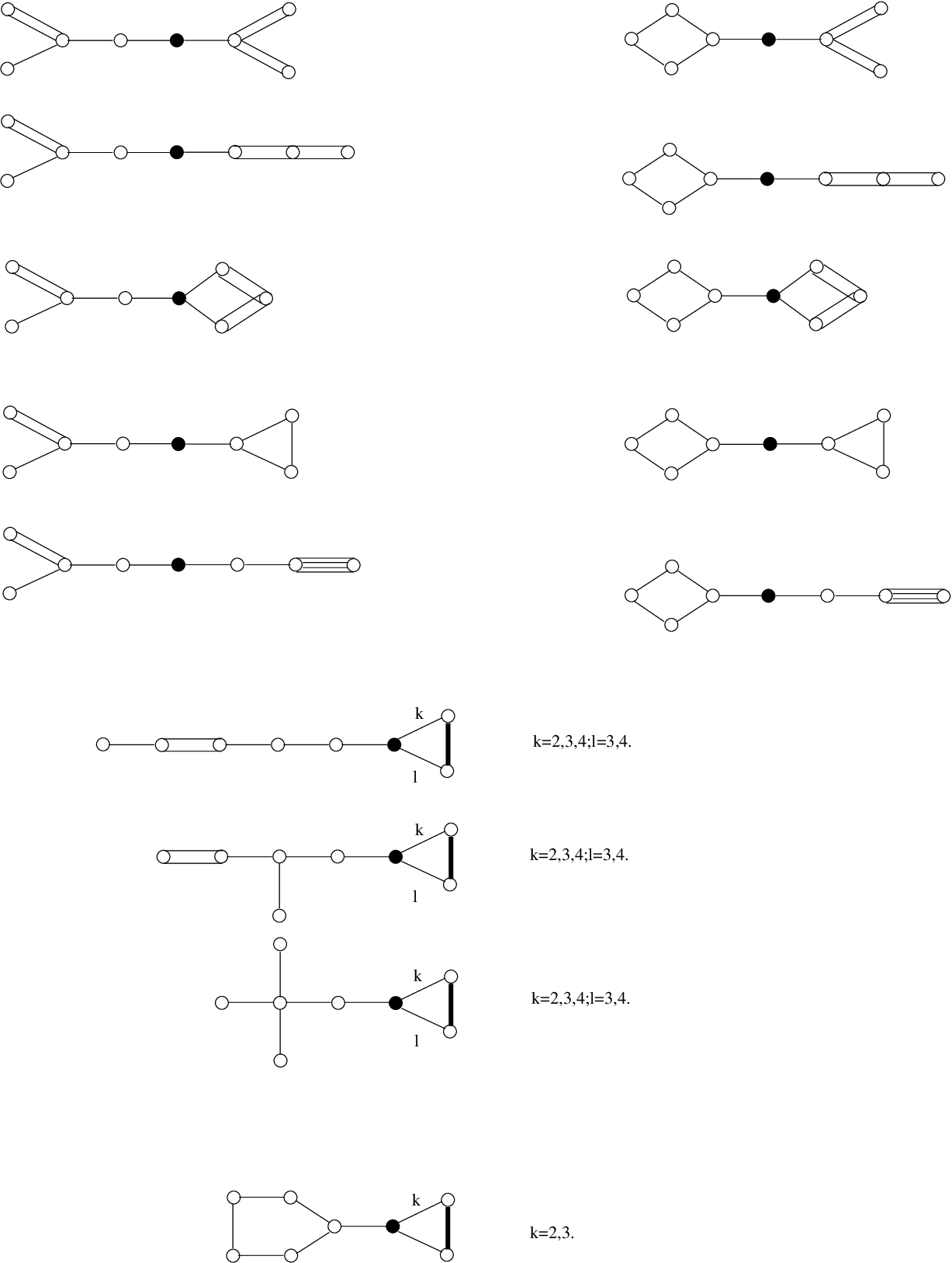,width=0.8\linewidth}
\end{center}
\caption{Pyramids in $\H^6$.}
\label{pyr6}
\end{table}

\begin{table}[htb!]
\begin{center}
\psfrag{k}{\scriptsize $k$}
\psfrag{k=2,3.}{\scriptsize $k=2,3.$}
\epsfig{file=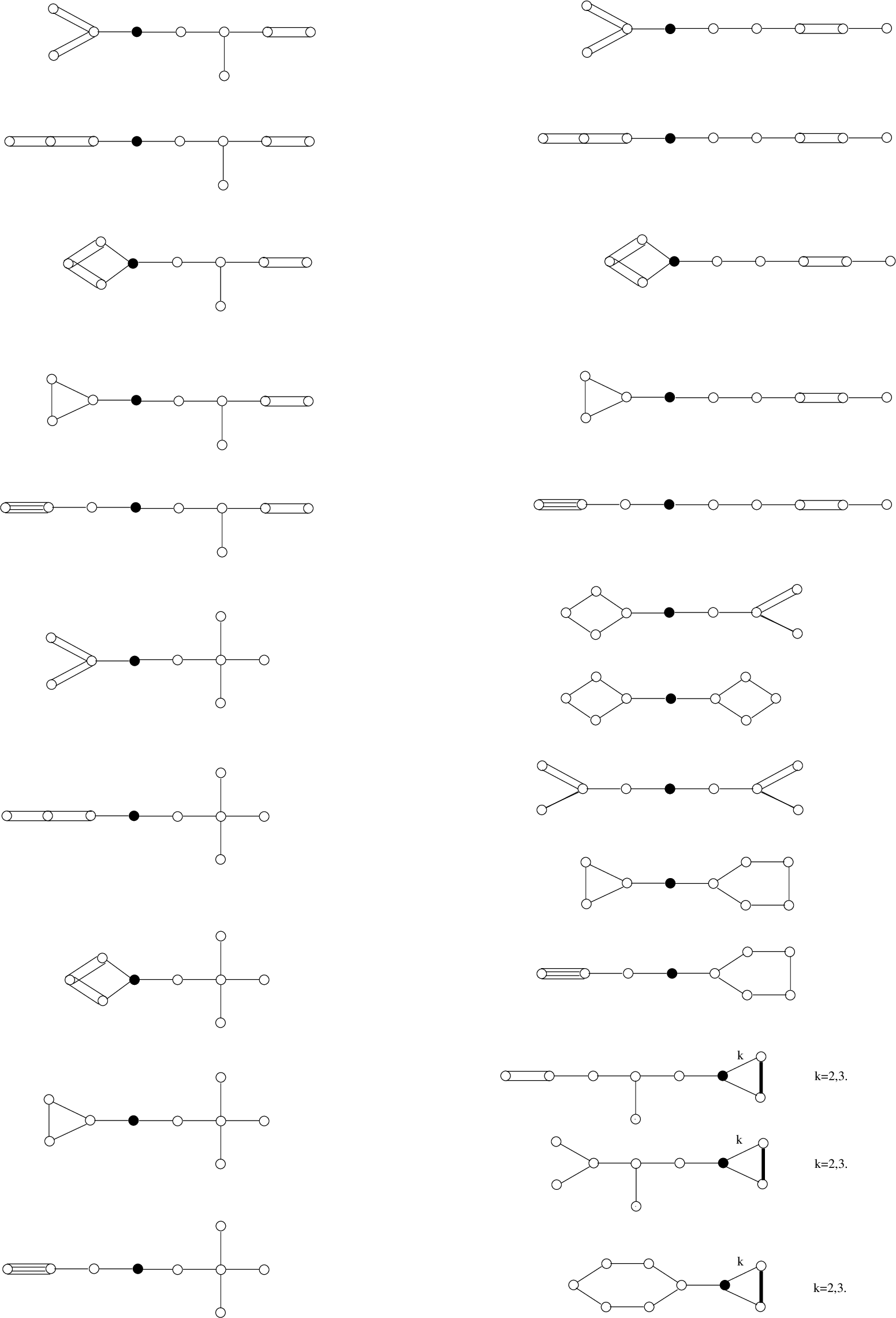,width=0.8\linewidth}
\end{center}
\caption{Pyramids in $\H^7$.}
\label{pyr7}
\end{table}

\begin{table}[htb!]
\begin{center}
\psfrag{k}{\scriptsize $k$}
\psfrag{k=2,3.}{\scriptsize $k=2,3.$}
\epsfig{file=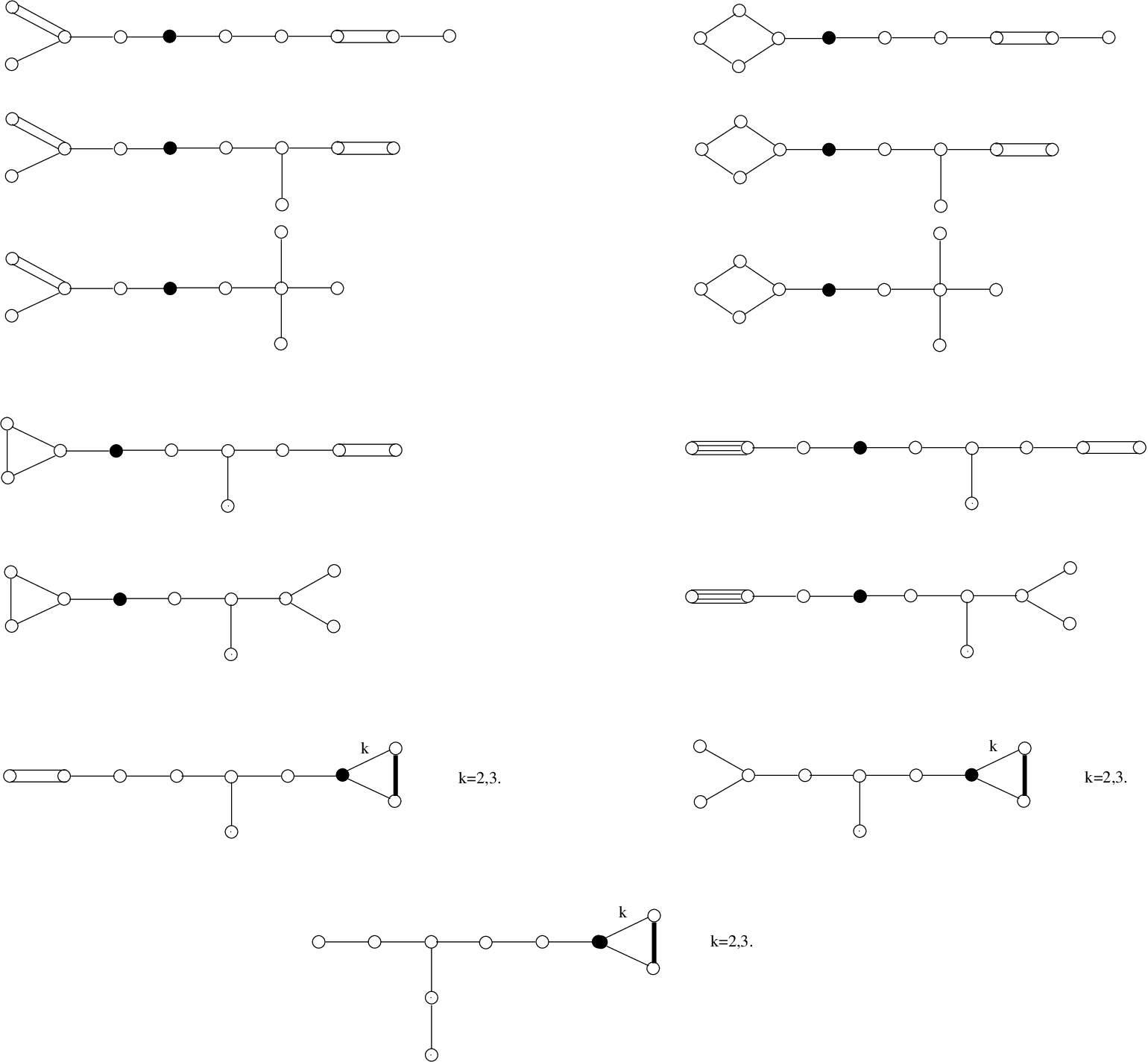,width=0.8\linewidth}
\end{center}
\caption{Pyramids in $\H^8$.}
\label{pyr8}
\end{table}

\begin{table}[htb!]
\begin{center}
\psfrag{k}{\scriptsize $k$}
\psfrag{k=2,3.}{\scriptsize $k=2,3.$}
\epsfig{file=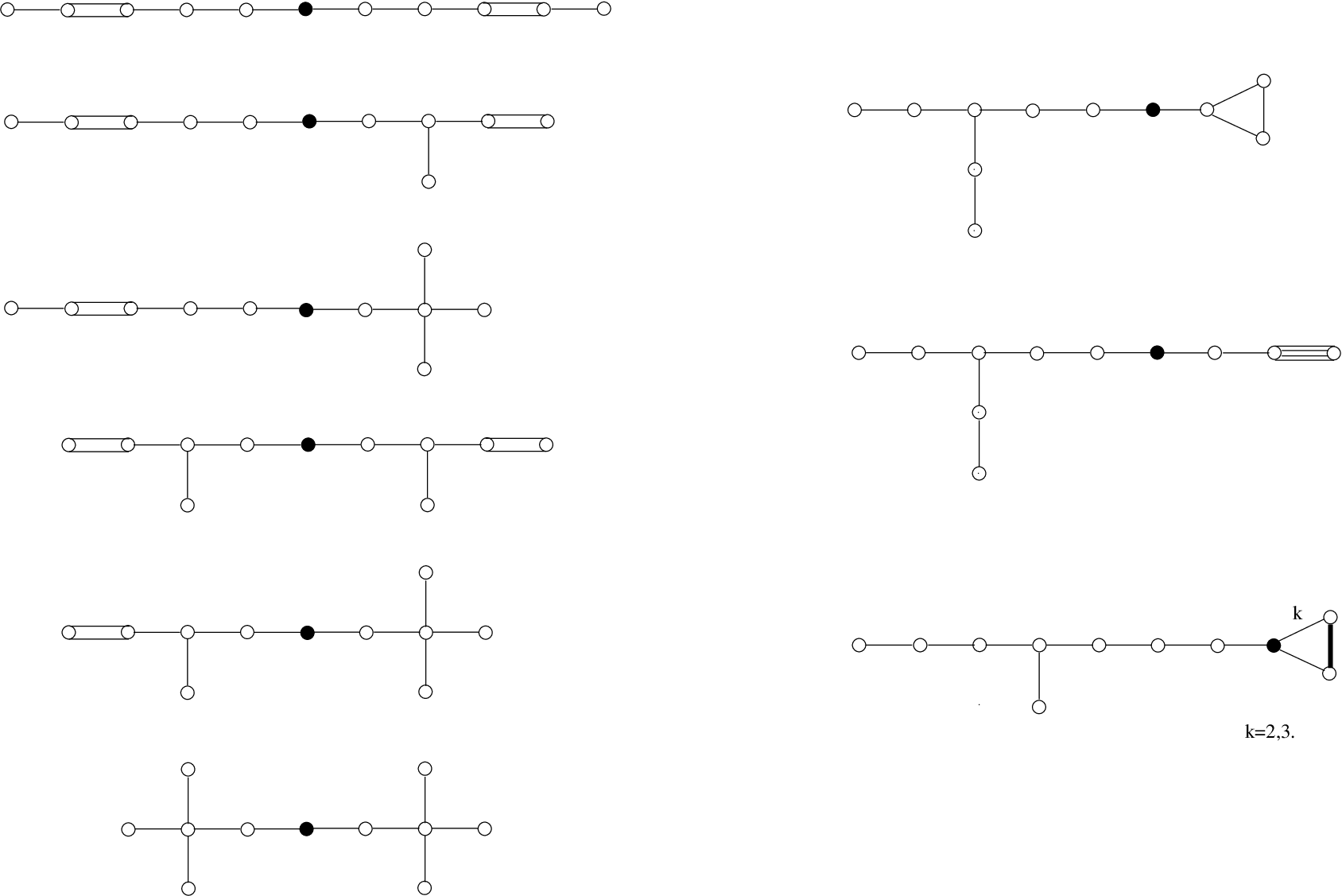,width=0.8\linewidth}
\end{center}
\caption{Pyramids in $\H^9$.}
\label{pyr9}
\end{table}

\begin{table}[htb!]
\begin{center}
\psfrag{k}{\scriptsize $k$}
\psfrag{l}{\scriptsize $l$}
\psfrag{l=3,4.}{\scriptsize $l=3,4.$}
\psfrag{k=2,3,4;}{\scriptsize $k=2,3,4;$}
\epsfig{file=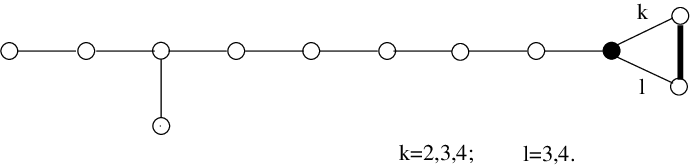,width=0.75\linewidth}
\end{center}
\caption{Pyramid in $\H^{10}$.}
\label{pyr10}
\end{table}

\begin{table}[htb!]
\begin{center}
\epsfig{file=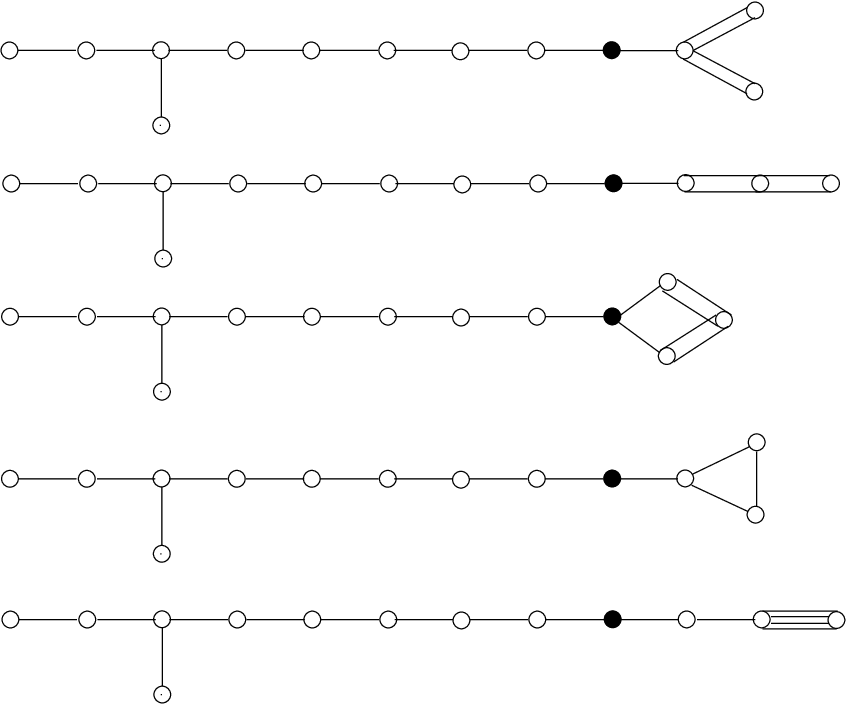,width=0.8\linewidth}
\end{center}
\caption{Pyramids in $\H^{11}$.}
\label{pyr11}
\end{table}

\begin{table}[htb!]
\begin{center}
\epsfig{file=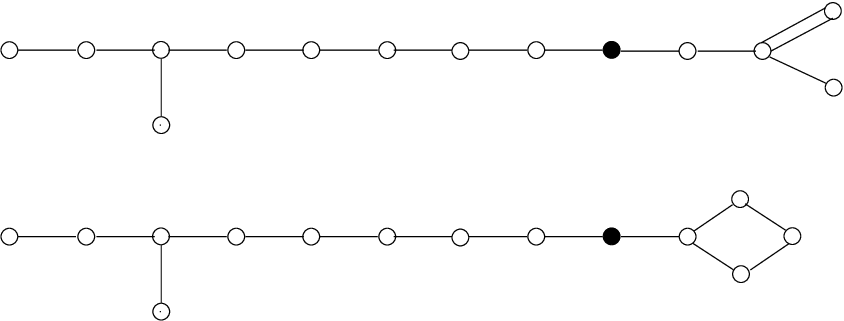,width=0.8\linewidth}
\end{center}
\caption{Pyramids in $\H^{12}$.}
\label{pyr12}
\end{table}

\begin{table}[htb!]
\begin{center}
\epsfig{file=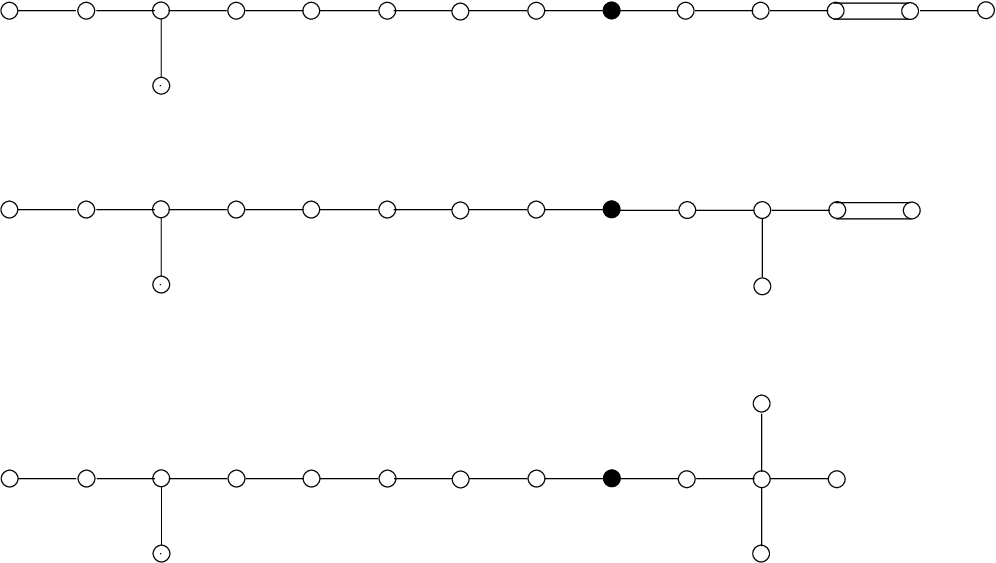,width=0.8\linewidth}
\end{center}
\caption{Pyramids in $\H^{13}$.}
\label{pyr13}
\end{table}

\begin{table}[htb!]
\begin{center}
\epsfig{file=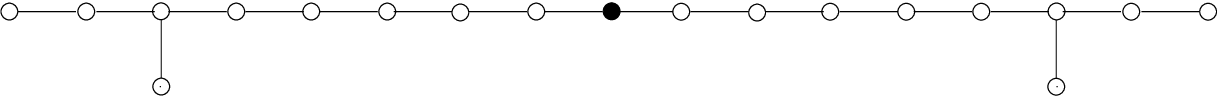,width=0.8\linewidth}
\end{center}
\caption{Pyramid in $\H^{17}$.}
\label{pyr17}
\end{table}

\end{document}